\documentclass[11pt,a4paper]{amsart}

\usepackage{amsmath}
\usepackage{amssymb}
\usepackage{latexsym}

\newcommand{\rpfeil}[2]{\stackrel{#2}{\verylongarrow{#1mm}}}

\newtheorem{theorem}{Theorem}[section]

\newtheorem{lemma}[theorem]{Lemma}
\newtheorem{proposition}[theorem]{Proposition}
\newtheorem{definition}[theorem]{Definition}
\newtheorem{corollary}[theorem]{Corollary}

\newtheorem{exmp}[theorem]{Example}
\newtheorem{exmps}[theorem]{Examples}
\newtheorem{rem}[theorem]{Remark}

\newcommand{\beeq}[1]{\begin{eqnarray}\label{#1}}
\newcommand{\eneq}{\end{eqnarray}}
\newcommand{\cal}{\mathcal} 
\newcommand  {\Br}      {\operatorname{Br}}
\newcommand{\ka}{{\cal A}}

\newcommand{\kc}{{\cal C}}

\newcommand{\ke}{{\cal E}}

\newcommand{\kk}{{\cal K}}

\newcommand{\km}{{\cal M}}
\newcommand{\ko}{{\cal O}}
\newcommand{\kp}{{\cal P}}

\newcommand{\kx}{{\cal X}}

\newcommand{\ZZ}{\mathbb{Z}}
\newcommand{\QQ}{\mathbb{Q}}
\newcommand{\RR}{\mathbb{R}}
\newcommand{\CC}{\mathbb{C}}

\newcommand{\PP}{\mathbb{P}}

\newcommand{\IC}{{\mathbb C}}

\newcommand{\IP}{{\mathbb P}}
\newcommand{\IQ}{{\mathbb Q}}
\newcommand{\IR}{{\mathbb R}}
\newcommand{\IZ}{{\mathbb Z}}

\renewcommand{\O}       {\mathcal{O}}

\newcommand{\verylongarrow}[1]{\hbox to #1{\rightarrowfill}}

\newcommand\mynote[1]{\marginpar{\ \\ \small \tt #1}}
\newcommand\bel[1]{{\mynote{#1}}\begin{equation}\label{#1}}
\newcommand\mylabel[1]{\label{#1}}

\begin{document}
\title[The Brauer group of K3 surfaces]{The Brauer group of  analytic K3 surfaces}

\author[Daniel Huybrechts]{Daniel Huybrechts}
\address{D.H.: Institut de Math\'ematiques de Jussieu, 2 place de Jussieu, 75251 Paris Cedex 05, France}
\email{huybrech@math.jussieu.fr}

\author[Stefan Schroer]{Stefan Schr\"oer}
\address{S.S.: Mathematische Fakult\"at, Ruhr-Universit\"at, 
         44780 Bochum, Germany}
\curraddr{Mathematisches Institut, Universit\"at K\"oln, Weyertal 86-90,
          50931 K\"oln, Germany}
\email{s.schroeer@ruhr-uni-bochum.de}

\begin{abstract}
We show that for complex analytic K3 surfaces any torsion class in 
$H^2(X,{\cal O}_X^*)$ comes from an Azumaya algebra. 
In other words, the Brauer group equals the cohomological Brauer group. 
For algebraic surfaces, such results go back to Grothendieck. 
In our situation, we use twistor spaces to deform a given 
analytic K3 surface to suitable projective K3 surfaces, and then 
stable bundles and hyperholomorphy conditions to pass back and 
forth between the members of the twistor family. 
\end{abstract}

\maketitle

In analogy to the isomorphism ${\rm Pic}(X)\cong H^1(X,{\cal O}_X^*)$,
Grothendieck investigated in \cite{GB II} the possibility of interpreting classes
in $H^2(X,{\cal O}_X^*)$ as geometric objects.
 He observed that the Brauer group $\Br(X)$, parameterizing equivalence classes of sheaves of Azumaya
 algebras on $X$, naturally injects into $H^2(X,\O_X^*)$. It is not difficult to see that 
 $\Br(X)\subset H^2(X,\O_X^*)$ is contained in the torsion part of $H^2(X,\O_X^*)$
 and Grothendieck asked:
 $${\it Is ~the~ natural ~injection ~\Br(X)\subset
 H^{\rm 2}(X,\O_X^*)_{\rm tor}~ an~ isomorphism?}$$

This question is of  interest in various geometric categories, e.g.\
$X$ might be a scheme, a complex space, a complex manifold, etc.  
It is also related to more recent developments in the application of 
complex algebraic
geometry to conformal field theory. Certain elements in $H^2(X,\O_X^*)$ have been
interpreted as so-called B-fields, and those are used to construct super conformal field
theories associated to Ricci-flat manifolds. Thus, understanding
the geometric meaning of the cohomological Brauer group $\Br'(X):= H^2(X,\O_X^*)_{\rm tor}$
is also of interest for the mathematical interpretation of string theory and mirror symmetry.

\bigskip

An affirmative answer to Grothendieck's question has been given only in very few special cases:

\medskip

$\bullet$ If $X$ is a complex curve, then $H^2(X,\O_X^*)=0$. Hence, $\Br(X)=\Br'(X)=H^2(X,\O_X^*)=0$
(see \cite[Cor.2.2]{GB II} for the general case of a curve).

$\bullet$ For smooth algebraic surfaces the surjectivity has been proved by Grothendieck
\cite[Cor.2.2]{GB II} and for normal algebraic  surfaces a proof was given
 more recently by Schr\"oer \cite{Schroeer}.

$\bullet$ Hoobler \cite{Hoobler} and Berkovich
\cite{Berkovitch} gave an affirmative answer for abelian varieties of any dimension
and Elencwajg and Narasimhan  gave another proof for complex tori
\cite{EN}.

$\bullet$ Demeyer and Ford showed equality of the Brauer group and the cohomological Brauer group
for smooth toric varieties \cite{Demeyer}.

$\bullet$ The case of smooth affine algebraic varieties 
was settled by Hoobler in \cite{Hoobler2}. Gabber proved the
surjectivity
in the case of arbitrary affine schemes and also for the union of two
affine schemes glued over an affine scheme \cite{Gabber}.
 
$\bullet$ Bogomolov and Landia 
\cite{Bogomolov; Landia 1990} showed that for a class $\alpha\in\Br'(X)$ 
on a scheme $X$, there is a blowing-up $f:Y\to X$ with 
$f^*(\alpha)\in\Br(Y)$. 

\medskip
Thus, besides the case of complex tori
the question has not been answered for any challenging class
of compact complex  varieties in dimension at least three. In fact, even in dimension
two, i.e.\ for surfaces, a complete answer is still missing for non-algebraic surfaces.
 
In this note we give a complete answer to  Grothendieck's question for 
analytic K3 surfaces.
\begin{theorem}\label{Main}
\mylabel{brauer k3}
Let $X$ be any
K3 surface. Then $\Br(X)=\Br'(X)=H^2(X,\O_X^*)_{\rm tor}$. 
\end{theorem}
 
If $0\leq\rho\leq 20$  is the Picard number of $X$ then one has
$\Br(X)=(\QQ/\ZZ)^{22-\rho}$.

For algebraic K3 surfaces this was known due to the aforementioned 
result of  Grothendieck.
The proof we give depends essentially on non-algebraic K3 surfaces and
their deformation theory.
Maybe the most interesting aspect of our approach is that the existence of the
Ricci-flat metric is crucial for the whole argument. 

Using standard techniques one easily shows that
our result implies an affirmative answer to Grothendieck's question
for all 
compact Ricci-flat K\"ahler surfaces (Cor. \ref{Ricci flat surfaces}). 

\section{Recollection: Brauer group}

In this first section we recall the basic facts concerning Brauer groups.
The standard reference for this is \cite{GB II}. As we will mainly
be interested in complex manifolds, we let $X$ be a connected complex manifold
endowed with the analytic topology. 

Let us begin with the definition of the Brauer group.

\begin{definition}
An \emph{Azumaya algebra} on the complex manifold $X$ is an associative (non-commutative) $\O_X$-algebra
${\cal A}$ which is locally (in the analytic topology)  isomorphic to a matrix algebra
${\rm M}_r(\O_X)$ for some $r>0$.
\end{definition}

Thus, any Azumaya algebra ${\cal A}$ is locally free of constant rank $r^2$. Two Azumaya algebras
are isomorphic if they are isomorphic as $\O_X$-algebras. By the Skolem-Noether theorem,
${\rm Aut}({\rm M}_r(\CC))\cong{\rm PGL}_r(\CC)$ (acting by conjugation). Hence, the set of isomorphism classes of
Azumaya algebras ${\cal A}$ of rank $r^2$ is in bijection with the set $H^1(X,{\rm PGL}_r(\O_X))$.

Note that any vector bundle $E$ of rank $r$ induces an Azumaya algebra ${\cal A}={\cal E}nd(E)$
of rank $r^2$. Moreover, the associated projective bundle $\PP(E)$ also defines a class
in $H^1(X,{\rm PGL}_r(\O_X))$ which actually coincides with the class defined by $\ka$.
Azumaya algebras arising this way will be declared `trivial' by  means of the following
equivalence relation.

\begin{definition}
Two Azumaya algebras ${\cal A}$ and ${\cal A}'$ are called \emph{equivalent} if there exist 
non-zero vector bundles
$E$ and $E'$ such that ${\cal A}\otimes{\cal E}nd(E)$ and ${\cal A}'\otimes{\cal E}nd(E')$ are isomorphic
Azumaya algebras.
\end{definition}

\begin{definition}
The \emph{Brauer group} $\Br(X)$ is the set of isomorphism classes of Azumaya algebras modulo the
above equivalence relation.
\end{definition}

That $\Br(X)$ is indeed a group stems from the fact
that ${\cal A}\otimes{\cal A}^{\rm op}\cong{\cal E}nd_{\O_X}(\cal A)$, i.e.\ ${\cal A}^{\rm op}$ defines an inverse
(with respect to "$\otimes$") of ${\cal A}$ as an element in $\Br(X)$. Clearly, any Azumaya
algebra of the form ${\cal A}\cong{\cal E}nd(E)$ defines  the unit element in $\Br(X)$.
 In fact, due to the following result, the converse holds also true.
\begin{lemma}
An Azumaya algebra $\ka$ of rank $r^2$ is trivial if and only if its class in
$H^1(X,{\rm PGL}_r(\O_X))$ is contained in the image of the natural map
$H^1(X,{\rm GL}_r(\O_X))\to H^1(X,{\rm PGL}_r(\O_X))$,
i.e.\ $\ka={\cal E}nd(E)$.
\end{lemma}

In order to prove this lemma one considers the long exact cohomology sequence 
induced by 
$$1\to\O_X^*\to{\rm GL}_r(\O_X)\to{\rm PGL}_r(\O_X)\to 1.$$
Since ${\rm GL}$ and ${\rm PGL}$ are not abelian, one only has a long exact sequence up to
$H^2(X,\ko_X^*)$: 

$$\to H^1(X,\O_X^*)\to H^1(X,{\rm GL}_r(\O_X))\to H^1(X,{\rm PGL}_r(\O_X))\to H^2(X,\O_X^* )$$
In particular, the kernel of the boundary map
$H^1(X,{\rm PGL}_r(\O_X))\to H^2(X,\O_X^*)$ is the set of isomorphism classes
of Azumaya algebras ${\cal E}nd(E)$ with $E$ locally free of rank $r$.
One then shows that all these maps factor over a group homomorphism
$$\delta:\Br(X)\to H^2(X,\O_X^*),$$
i.e.\ for any $r$ one has a commutative diagram
$$\begin{array}{ccc}
H^1(X,{\rm PGL}_r(\O_X))&\to&H^2(X,\O_X^*)\\
\downarrow&&\parallel\\
\Br(X)&\to&H^2(X,\O_X^*)\\
\end{array}$$
Together with the exact sequence above  one finds that
$\delta:\Br(X)\to H^2(X,\O_X^*)$ is injective.

In order to see that the image of $\delta$ is contained in the torsion part one uses
the commutative diagram
$$\begin{array}{ccccccccc}
1&\to&\mu_r&\to&{\rm SL}_r&\to&{\rm PGL}_r&\to&1\\
&&\cap&&\cap&&\parallel&&\\
1&\to&\,\,\O_X^*&\to&{\rm GL}_r&\to&{\rm PGL}_r&\to&1
\end{array}$$
The cohomology sequence induced by the short exact sequence on the bottom
provides us with a boundary map
$\eta_r:H^1(X,{\rm PGL}_r(\O_X))\to H^2(X,\mu_r)$.
Since the map $H^1(X,{\rm PGL}_r(\O_X))\to H^2(X,\O_X^*)$
factorizes over $\eta_r$ the image of it is $r$-torsion.
Thus, one obtains the natural injection
$$\delta:\Br(X)\hookrightarrow H^2(X,\O_X^*)_{\rm tor}$$ and
$H^2(X,\O_X^*)_{\rm tor}$ is called the \emph{cohomological Brauer group} $\Br'(X)$.

Later we will also make use of the Kummer sequence
$$1\to\mu_r\to\O^*_X\rpfeil{5}{(~)^r}\O_X^*\to 1$$
and the induced short exact sequence
\begin{equation}
1\to H^1(X,\O_X^*)/r\cdot H^1(X,\O_X^*)\to H^2(X,\mu_r)\to H^2(X,\O_X^*)_{r{\rm -tor}}\to1.
\label{Kummer}
\end{equation}
Passing to the direct limits we obtain
a surjection
$$ H^2(X,\QQ/\ZZ)\rpfeil{5}{\exp}\lim H^2(X,\mu_r)\twoheadrightarrow H^2(X,\O_X^*)_{\rm tor}.$$
The sheaf $\mu_r$ viewed as $\ZZ/r\ZZ$ also sits in the exact sequence
$$0\to\ZZ\rpfeil{5}{{\cdot r}}\ZZ\to\mu_r\to 0.$$
The induced boundary morphism $\beta:H^2(X,\mu_r)\to H^3(X,\IZ)$
is called the Bockstein.
The Bockstein is compatible with the boundary map
$H^2(X,\O_X^*)\to H^3(X,\ZZ)$ induced by the exponential sequence.
This is due to the commutative diagram
$$\begin{array}{ccccccccc}
0&\to&\ZZ&\rpfeil{5}{{\cdot r}}&\ZZ&\to&\mu_r&\to&0\\
&&\parallel&&\downarrow\cdot {1}/{r}&&\cap&&\\
0&\to&\ZZ&\to&\O_X&\rpfeil{5}{\exp}&\O_X^*&\to&0\\
\end{array}$$

Thus, a class in $\alpha\in H^2(X,\mu_r)$ induces a topologically trivial Brauer class,
i.e. its image in $H^3(X,\IZ)$ is trivial,
if and only if its Bockstein $\beta(\alpha)$ is trivial.

We will also need the following easy 

\begin{lemma}\mylabel{det versus brauer}
Let $E$ be a vector bundle of rank $r$ on the complex manifold $X$. 
Then the image of $\det(E)^*$ under the natural mapping $H^1(X,\O_X^*)\to H^2(X,\mu_r)$ 
equals $\eta_r(\IP(E))$.
\end{lemma}

\proof Indeed, if $E$ is given by the cocycle $\{\varphi_{ij}\in\Gamma(U_{ij},{\rm GL}_r)\}$
such that there exist $\lambda_{ij}\in\Gamma(U_{ij},\O^*_X)$
with $\lambda_{ij}^r=\det(\varphi_{ij})$, then $\det(E)^*\in H^1(X,\O_X^*)$ is given
by $\{\det(\varphi_{ij})^{-1}\}$ and its image in $H^2(X,\mu_r)$ by $\{\lambda^{-1}_{ij}\lambda_{jk}^{-1}
\lambda_{ik}\}$. 

On the other hand, $[\PP(E)]=\{\overline{\varphi_{ij}}\}
=\{\overline{\varphi_{ij}\cdot\lambda_{ij}^{-1}}\}\in H^1(X,{\rm PGL}_r(\O_X))$
and $\varphi_{ij}\cdot\lambda^{-1}_{ij}\in{\rm SL}_r$. This shows that the image of
$[\PP(E)]$ in $H^2(X,\mu_r)$ is given by the cocycle $\{\varphi_{ij}\lambda_{ij}^{-1}
\varphi_{jk}\lambda_{jk}^{-1}(\varphi_{ik}\lambda_{ik}^{-1})^{-1}\}$.\qed

\bigskip
In case that $X$ is projective one may as well work with the \'etale
topology. The Kummer sequence immediately shows that the two cohomological Brauer groups
coincide.
Eventually, note that any Brauer class on a K3 surface is automatically
topologically trivial, for $H^3(X,\IZ)=0$.
\section{Recollection: K3 surfaces}

As a reference for the theory of K3 surfaces we recommend \cite{Periodes}.
By definition a K3 surface is a compact complex surface $X$ such that the canonical
bundle $K_X$ is trivial and $H^1(X,\ko_X)=0$. It has been shown by Siu that any K3 surface
is K\"ahler. Thus, the K\"ahler cone $\kk_X\subset H^{1,1}(X,\IR)$ of all
K\"ahler classes is a non-empty open convex  cone. Moreover, due to Yau's solution of the
Calabi-conjecture any class in $\omega\in\kk_X$ can uniquely be represented by a Ricci-flat K\"ahler
form, which we will also call $\omega$.

For holomorphic symplectic (also called hyperk\"ahler) manifolds in general and for
K3 surfaces in particular any Ricci-flat K\"ahler metric $g$ is actually hyperk\"ahler.
More precisely this means that in addition to the given complex structure $I$  defining $X$
there exist complex structures $J$ and $K=I\circ J=-J\circ I$ all making $g$ a K\"ahler metric.
The induced sphere of complex structures $\{aI+bJ+cK~|~a^2+b^2+c^2=1\}$ is isomorphic 
to a projective line called $\IP(\omega)\cong\IP^1$. 

Twistor theory of K3 surfaces shows that
there exists a universal complex structure on $X\times \IP(\omega)$. This complex manifold
will be called $\kx(\omega)$. Moreover, it is known that the projection
$\kx(\omega)\to\IP(\omega)$ is holomorphic, and by definition the fiber over a complex
structure $\lambda$ is isomorphic to $X$ endowed with this structure. 
Usually we will denote by $0\in\IP(\omega)$ the point that corresponds to the original
complex structure.

Although the twistor space $\kx$ itself is not k\"ahler, each fiber $\kx_t$ comes
with a natural K\"ahler form $\omega_t$. If $t$ corresponds to the complex structure
$\lambda_t$ then $\omega_t=g(\lambda_t(~),~)$. In particular, $\omega_0=\omega$.

\begin{lemma}\mylabel{11 on Xt}
Let $\alpha\in H^2(X,\IZ)$ with $\alpha^2>0$. Then there exists $t\in \IP(\omega)$ such that
$\alpha$ is of type $(1,1)$ on $\kx_t$ and $\kx_t$ is projective.
\end{lemma}

\proof We denote the holomorphic two-form (or rather its cohomology class)
on $\kx_t$ by $\sigma_t$. Thus,
$Q_\omega:=\{\sigma_t\}\subset \IP(\IC\omega_I\oplus\IC\omega_J\oplus \IC\omega_K)$
is the quadric defined by the intersection form restricted to the three-dimensional
space spanned by $\omega_I,\omega_J,\omega_K$.
The class $\alpha$ is of type $(1,1)$ with respect to the complex structure
corresponding to $t\in\IP(\omega)$ if and only if $\alpha\wedge\sigma_t=0$.
Since the intersection of the quadric $Q_\omega$ with the hyperplane
defined by $\alpha \wedge x=0$ is non-empty, one finds a point $t$ as required.

The second assertion follows from the general fact that a compact complex
surface is projective if and only if there exists a line bundle
$L$ with ${\rm c}_1(L)^2>0$. This is applied to the line bundle $L$ whose
first Chern class is $\alpha$.
\qed

\bigskip

We next discuss stable and hyperholomorphic bundles. Stability will always mean
slope-stability and if we write $\omega_t$-stable it means slope-stable with respect
to the K\"ahler class $\omega_t$ on $\kx_t$.

\begin{definition}
A vector bundle $E$ (or, projective bundle $P$) on $\kx_t$ is called \emph{hyperholomorphic}
(with respect to the twistor space $\kx(\omega)$) if there exists a vector bundle $\ke$ (respectively,
projective bundle $\kp$) on $\kx(\omega)$ such that $\ke|_{\kx_t}\cong E$ (respectively,
$\kp|_{\kx_t}\cong P$).
\end{definition}

Since the twistor space $\kx(\omega)$ comes with the natural $\kc^\infty$-trivialization
$\kx(\omega)\cong X\times\IP(\omega)$, any holomorphic vector bundle
$E$ on a fiber $\kx_t$ lives naturally as a complex (!) vector bundle on any
other fiber $\kx_s$. Thus, $E$ is hyperholomorphic if and only if the complex vector bundle
$E$ admits a $\bar\partial$-operator on any $\kx_s$ depending holomorphically on $s$.

The following result is due to Verbitsky. For the convenience of the reader we include
a sketch of its proof.

\begin{proposition}\mylabel{Verbitsky Prop}
Let $E$ be an $\omega_t$-stable vector bundle on $\kx_t$. Then
the associated projective bundle $\IP(E)$ is hyperholomorphic with respect
to $\kx(\omega)$ (see \cite{Verbitsky}).
\end{proposition}

\proof
To give an idea of the proof we first assume that the determinant of $E$
is trivial. Then we show that the vector bundle $E$  itself is hyperholomorphic
and hence also $\IP(E)$. Due to Donaldson's result on the existence of Hermite-Einstein
metrics there exists a hermitian metric on $E$ such that the curvature
$F_\nabla\in\ka^{2}(X,{\rm End}(E))$ of the induced Chern connection
$\nabla$ on $E$ satisfies $\Lambda_{\omega_t}F_\nabla=0$ (this is the Hermite-Einstein condition).
Since we are working on a surface, this condition is equivalent to $F\wedge \omega=0$.
We are going to show that the connection $\nabla$ is hyperholomorphic, i.e.\
with respect to any other complex structure given by a point $s\in\IP(\omega)$
the $(0,1)_s$-part $\nabla^{(0,1)_s}$ of $\nabla$ is a $\bar\partial$-operator for the complex bundle
$E$ on $\kx_s$. Indeed, $(\nabla^{(0,1)_s})^2$ is the $(0,2)_s$-part of the curvature
$F_\nabla$. By definition of the Chern connection, $F_\nabla$ is of type $(1,1)$ on
$\kx_t$. Thus, its suffices to show that it stays of type $(1,1)$ when $\kx_t$ is deformed to
$\kx_s$.

In fact, it suffices to show that $F$ is of type $(1,1)$ for $J$ and $K$, because
then it will automatically be of type $(1,1)$ with respect to any complex structure parameterized
by $\IP(\omega)$.
Since the curvature $F$ is real, it is of type $(1,1)$ with respect to $J$ if and only if
$F\wedge\sigma_J=0$, where $\sigma_J$ is the holomorphic two-form
with respect to the complex structure $J$. Using $\sigma_J=\omega_K+i\omega_I$ the Hermite-Einstein
condition immediately yields $F\wedge{\rm Im}(\sigma_J)=0$.
On the other hand, ${\rm Re}(\sigma_J)=\omega_K={\rm Im}(\sigma_I)$. Hence,
$F\wedge{\rm Re}(\sigma_J)=F\wedge{\rm Im}(\sigma_I)=0$, since $F$ is of type
$(1,1)$ with respect to $I$. Analogously one proves $F\wedge\sigma_K=0$.

If the determinant of $E$ is no longer trivial, then the Hermite-Einstein condition
for $E$ reads $F\wedge \omega_t=\omega^2\cdot\mu\cdot{\rm id}_E$, where
$\mu$ is a constant measuring the degree of $\det(E)$.
The given Chern connection $\nabla$ for $E$ on $\kx_t$ defines a holomorphic structure
on the induced projective bundle $\IP(E)$ with respect to another
complex structure $s\in\IP(\omega)$ if and only if the tracefree part $F_0$ of
$F$, given by $F=F_0+\frac{{\rm tr}(F)}{r}{\rm id}_E$,  is of type $(1,1)$ with respect to $s$.
In order to see this we replace $F$ by $F_0$ in the argument above and use the Hermite-Einstein
condition.
\qed

\bigskip

\section{Proof}

The rough idea of the proof of Theorem \ref{Main} is as follows.
In order to show the assertion it suffices to realize any element
$\bar\alpha\in H^2(X,\mu_r)$ as the class $\eta_r(\ka)$ of some Azumaya algebra.
Here we use the surjectivity of the map 
$H^2(X,\mu_r)\to H^2(X,\ko_X^*)_{r{\rm -tor}}$.
Clearly, if $\ka_t$ is a flat family of Azumaya algebras over a deformation
$\kx_t$ of $\kx_0=X$, then $\eta_r(\ka_t)$ is constant. Thus, we may try to first deform
$X$ to some $\kx_t$ and prove the existence of the required Azumaya algebra there.
In order to be able to deform back to $X$ we have to find a `hyperholomorphic'
Azumaya algebra on $\kx_t$. The reason why it might be easier to find an Azumaya algebra
on some deformation is the fact that the Picard group changes when $X$ is deformed
and one can arrange things such that the given class in $ H^2(X,\mu_r)$ is actually in 
the image of $H^1(X,\ko_X^*)\to H^2(X,\mu_r)$ induced by the Kummer sequence.
Thus, the complex structure is deformed in a way that the Brauer class induced
by $\bar\alpha$ becomes trivial. It is a curious fact that by making the assertion trivial
on some deformation of $X$ one can in the end conclude the existence of a
non-trivial Azumaya algebra on $X$ itself.
Note that in his article \cite{dJ}, de Jong  also used the fact that under deformation of the surface
a trivial Azumaya algebra sometimes becomes non-trivial.

\medskip

Let us now pass on to the details of the proof. 
Let $\beta\in \Br'(X)$ be $r$-torsion and choose a lift $\bar\alpha\in H^2(X,\mu_r)$. 

\begin{lemma}
There exists an element $\alpha\in H^2(X,\IZ)$ inducing $\bar\alpha$ such that
$\alpha^2>0$.
\end{lemma}

\proof
This is elementary: Choose any class $\alpha$  that lifts $\bar\alpha$ and add
$k\cdot r\cdot \beta$ with $k\gg0$ and $\beta\in H^2(X,\IZ)$ an arbitrary class
with $\beta^2>0$. 
\qed

\bigskip
 
 The next lemma is a refinement of Lemma \ref{11 on Xt}
\begin{lemma}\mylabel{defo picard one}
If $\omega\in\kk_X$ is a very general K\"ahler class on $X$ then there exists a point 
$t$ in the twistor base $\IP(\omega)$ such that the fiber $\kx_t$ is projective with
Picard number one
and $\alpha$ is of type $(1,1)$ on $\kx_t$.
\end{lemma}

\proof This is proved by standard techniques. We nevertheless provide
all details.
Let $V_\RR=(\kk_X-\left\{0\right\})/\RR^*$ 
be the  real projectivization of the K\"ahler cone
$\kk_X$. It is an open subset of the real projective space 
$(H^{1,1}(X,\RR)-\left\{0\right\})/\RR^*$ and is contained in the
complex manifold $V=(\CC\kk_X-\left\{0\right\})/\IC^*$,
the complex projectivization of the open convex cone
$\CC\kk_X=\kk_X+\sqrt{-1}H^{1,1}(X,\RR)\subset H^{1,1}(X,\CC)$. 
Locally $V_\RR\subset V$  looks like 
$\RR^{20}\subset\CC^{20}$. 

Fix some K\"ahler class $[\omega_0]\in  \kk_X$.
We now construct
an open neighborhood $U\subset V$ of $[\omega_0]$ 
and an open embedding $\varphi:U\to S_\alpha$, where
$S_\alpha$ is the hypersurface in the moduli space $\km$ of marked K3 surfaces
$(Y,\psi:H^2(Y,\IZ)\cong H^2(X,\IZ))$ such that $\psi^{-1}(\alpha)$ is of type $(1,1)$ on $Y$. 

By definition $S_\alpha$ is a hyperplane section of the moduli space $\km$ and thus
intersects any `generalized twistor line'
$\IP(\omega):=\kp^{-1}(\IP(\IC\omega\oplus\IC\sigma\oplus\IC\bar\sigma))$ for any $\omega\in V$
close to $\omega_0$ either in a non-reduced point or in two different points.
Here, $\kp$ is the period map $\kp  :\km\to\IP( H^2(X,\IC))$.
Let $D\subset V$ be the subset of all $[\omega]$ such that 
$\PP(\omega)\cap S_\beta$ is non-reduced. This is a closed complex
subspace which does not intersect $V_\IR$. 
Thus, we may find an open neighborhood $U\subset V\setminus D$ of $[\omega_0]$.
Choosing locally one of the two intersection points yields an open embedding 
$\varphi:U\to S_\beta$.

Now set $S:=\cup S_{\beta}$, where the union, which is countable, runs over all classes
$\beta\in H^2(X,\ZZ)$ that
are $\QQ$-linear independent from $\alpha\in H^2(X,\ZZ)$.
Set $u=[\omega_0]\in U$ and $s=\varphi(u)\in S_\beta$.
Suppose that $\varphi(U_\RR)\subset S$. It then follows for the real tangent
spaces that
$$
\varphi_*(T_u U_\RR)\subset
\bigcup T_s(S_\alpha\cap S_{\beta})
$$
and hence $\varphi_*(T_u U_\RR)\subset
T_s(S_\alpha\cap S_{\beta})$ for some $\beta$. 
Clearly, $T_s(S_\alpha\cap S_{\beta})\subset T_sS_\alpha$ is a 
complex subspace of complex codimension $\geq 1$, and moreover
the induced map $\varphi_*:T_uU=T_u U_\RR\otimes\CC\to T_sS_\alpha$ factors over 
$T_s(S_\alpha\cap S_{\beta})$. On the other hand,
the $\IC$-linear map $\varphi_*:T_uU\to T_sS_\alpha$ 
is bijective, and this  yields a contradiction.

Consequently,  for the very general  $[\omega]\in U_\RR$ the image 
$\varphi([\omega])\in S_\alpha$ is not contained in
any $S_{\beta}$. Hence, $\omega$ and $t:=\varphi([\omega])$ satisfy the assertion.
\qed

\bigskip

In particular the K3 surface $\kx_t$ as in the previous lemma will be
projective and $\alpha={\rm c}_1(L)$ for some holomorphic
line bundle $L$. Furthermore, $\kx_t$ is endowed with
a natural K\"ahler structure $\omega_t$.

\begin{proposition}
There exists an $\omega_t$-stable vector bundle $E$ on $\kx_t$ of rank $r$ and 
with $\det(E)=L^*$.
\end{proposition}

\proof
Since the Picard number of $\kx_t$ is one, there exists only one polarization $H$
(up to scaling). It is a standard result (eg.\ \cite{HL}) that for a given
polarization $H$ and determinant there exists a $H$-stable vector bundle
(for high enough second Chern number). Since the Picard number is one,
a $H$-stable vector bundle will in fact be stable with respect to any K\"ahler class.
This way we find plenty of $\omega_t$-stable vector bundles $E$ with
$\det(E)\cong L^*$.
\qed

\bigskip

{\bf End of the proof of Theorem \ref{brauer k3}.}
Due to Proposition \ref{Verbitsky Prop} we have found a holomorphic projective bundle
$\kp$ on $\kx(\omega)$ such that $\kp|_{\kx_t}\cong\IP(E)$, where
$E$ is a holomorphic vector bundle on $\kx_t$ of rank $r$ and with ${\rm c}_1(E)=-\alpha$.

Due to Lemma \ref{det versus brauer} the class $\eta_r(\IP(E))$ equals $\bar\alpha$.
Hence, for any $s\in\IP(\omega)$ one has $\eta_r(\kp_s)=\bar\alpha$.
In particular, $\bar\alpha$ is realized by a projective bundle and hence by an
Azumaya algebra on $X=\kx_0$.

This eventually yields that the given cohomological Brauer
class $\beta$ on $X$ is the class of an Azumaya algebra on $X$,
and the proof of Theorem \ref{brauer k3} is complete.

\medskip
In order to see that $\Br(X)\cong (\IQ/\IZ)^{22-\rho}$ use
the exponential sequence. The kernel of 
$H^2(X,\IZ)\to H^2(X,\ko_X)\cong\IC$ is the
Picard group and $H^2(X,\IZ)$ is of rank $22$. Hence, the torsion
part of $H^2(X,\ko_X)/H^2(X,\IZ)$ is isomorphic to $(\IQ/\IZ)^{22-\rho}$. 
Actually, the latter argument has
not much to do with K3 surfaces, and the same reasoning
shows that
the topologically trivial part of $\Br'(X)$ equals $(\QQ/\ZZ)^{b_2-\rho}$
for any compact  complex space.

\section{Further comments} 

{\bf 1.} The result for K3 surfaces is enough to prove the following 

\begin{corollary}\label{Ricci flat surfaces}
Let $X$ be a Ricci-flat compact K\"ahler surface, i.e.\ $0={\rm c}_1(X)\in H^2(X,\IR)$.
Then $\Br(X)=\Br'(X)$.
\end{corollary}

\proof In fact any compact Ricci-flat K\"ahler manifold $X$ admits a 
finite \'etale cover $\tilde X\to X$ such that $\tilde X$ is isomorphic to
a product of complex tori, hyperk\"ahler, and Calabi-Yau manifolds (cf.\ \cite{Beauville}).
Thus, our surface $X$ admits a finite \'etale cover $\pi:\tilde X\to X$ such that
$\tilde X$ is either a K3 surface or a complex torus.
For complex tori due to \cite{Hoobler, Berkovitch, EN} and for K3 surfaces
 due to Thm.\ \ref{brauer k3}
one knows that $\Br(\tilde X)\cong \Br'(\tilde X)$. Using a result of Gabber for finite maps
\cite[Ch.II, Lemma 4]{Gabber} one concludes from this $\Br(X)=\Br'(X)$.\qed

\bigskip

{\bf 2.} It is tempting to try to extend the above proof to the case of compact hyperk\"ahler manifolds,
which generalize the notion of K3 surfaces in many respects. In fact, most of the arguments go
through. E.g.\ Lemma \ref{defo picard one} can be proved in arbitrary dimension by using results of
\cite{H}, in particular the projectivity criterion.
One could hope to prove this way that any topologically trivial Brauer class
on a compact Ricci-flat K\"ahler manifold is contained in the Brauer group.
The point were the proof as presented is insufficient for higher-dimensional hyperk\"ahler manifolds
is in Proposition \ref{Verbitsky Prop}. 
In the higher-dimensional version of this result, also due to Verbitsky,
one needs to assume that the stable vector bundle $E$ with first Chern class
$\alpha$ has a discriminant $2r{\rm c}_2(E)-(r-1){\rm c}_1(E)^2$ of type $(2,2)$ with respect
to any complex structure parametrized by the twistor space in question. Of course, in complex
dimension two, the latter condition is automatic, since $H^{2,2}(X)=H^4(X)$.
In higher dimension however, $H^{2,2}(X)$ is strictly smaller than $H^4(X)$ and, moreover,
depends on the complex structure. To show the existence of the required stable bundle seems highly non-trivial.
In fact, the only hyperholomorphic projective bundles on higher-dimensional hyperk\"ahler
manifolds known to us are all descendants of the tangent bundle (which is stable!) and
we were not able to produce a bundle that would serve our purpose.

As before, the result for topologically trivial Brauer classes on hyperk\"ahler mani\-folds together
with the result of Hoobler, Berkovich, Elencwajg, and Narasimhan would show that
any topologically trivial Brauer class on a compact Ricci-flat K\"ahler variety is
contained in the Brauer group. (One uses that by definition $H^2(X,\ko_X)=0$ for a true Calabi-Yau manifold.)
This would certainly be strong  evidence for an affirmative
answer of Grothendieck's question.
Note however that topologically non-trivial classes on Calabi-Yau threefolds are
of particular interest in M-theory and for those the question remains open.


{\footnotesize \begin{thebibliography}{mm}

\bibitem{Periodes} \em G\'eom\'etrie des surfaces K3: modules et p\'eriodes.
\em S\'eminaires Palaiseau. ed A.\ Beauville, J.-P.\ Bourguignon, M.\ Demazure.
Ast\'erisque 126 (1985).

\bibitem{Beauville}
A.\ Beauville \em Vari\'et\'es K\"ahleriennes dont la premi\`ere classe de Chern est nulle. \em
J.\ Diff.\ Geom.\ 18 (1983) 755-782.

\bibitem{Berkovitch}
V.\ G.\ Berkovitch \em 
The Brauer group of abelian varieties. \em
Funkcional Anal. i Prilo\v zen. 6 (1972), 10-15.
 
\bibitem{Bogomolov; Landia 1990}
F.~Bogomolov, A.~Landia \em
$2$-cocycles and Azumaya algebras under birational transformations of 
algebraic schemes. \em
Compositio Math.\ 76 (1990),  1-5. 

\bibitem{Demeyer}
F.\ Demeyer, T.\ Ford
\em  On the Brauer group of toric varieties. \em
Trans.\ Amer.\ Math.\ Soc.\  335  (1993), 559-577. 

\bibitem{EN}
G.\ Elencwajg, S.\ Narasimhan \em
Projective bundles on a complex torus. \em 
J.\ Reine Angew.\ Math.\ 340 (1983), 1-5.

\bibitem{Gabber}
O.\ Gabber \em Some theorems on Azumaya algebras. \em 
In: M.~Kervaire, M.~Ojanguren (eds.), Groupe de Brauer, pp. 129--209,
Lecture Notes in Math.\ 844.
Springer, Berlin (1981). 

\bibitem{GB II} A.\ Grothendieck \em
Le groupe de Brauer. \em
In: J.~Giraud (ed) et al.: Dix expos\'es sur la cohomologie des
sch\'emas, pp.\ 46-189.
North-Holland, Amsterdam (1968). 


\bibitem{Hoobler}
R.\ Hoobler \em Brauer groups of abelian schemes. \em Ann.\ Sci. \'Ecole Norm.\ Sup.\ 5  (1972),
45-70. 

\bibitem{Hoobler2}
R.\ Hoobler \em  A cohomological interpretation of Brauer groups of rings.  \em
Pacific J.\ Math.\  86  (1980),  89-92.

\bibitem{HL}
D.\ Huybrechts, M.\ Lehn \em
The geometry of moduli spaces of sheaves. \em 
Aspects of Math.\ E31. 
Vieweg, Braunschweig, (1997).

\bibitem{H}
D.\ Huybrechts
\em Compact hyper-K\"ahler manifolds: basic results. \em
Invent.\ Math.\  135  (1999),  63-113.
Erratum: Invent.\ math.\ 152 (2003), 209-212

\bibitem{dJ}
J.\ de Jong
\em The period-index problem for the brauer group of an
algebraic surface. \em
Preprint. {\tt http://www-math.mit.edu/$\sim$dejong/}

\bibitem{Schroeer}
S.\ Schr\"oer \em There are enough Azumaya algebras on surfaces.  \em
Math.\ Ann.\ 321 (2001), 439-454.

\bibitem{Siu}
Y.-T.\ Siu \em Every K3 surface is K\"ahler. \em Invent.\ math. 73 (1983), 139-150.

\bibitem{Verbitsky}
M.\ Verbitsky \em Hyperholomorphic bundles over a hyper-K\"ahler manifold. \em
J.\ Alg.\  Geom.\  5  (1996),  633-669. 

\bibitem{Yau} S.-T.\ Yau \em On the Ricci curvature of a compact K\"ahler manifold
and the complex Monge-Amp\`ere equation. \em Comm.\ Pure Applied Math.\ 31 (1978), 339-411.

\end{thebibliography}}
\end{document}